\documentclass[english]{extarticle}
\usepackage[T1]{fontenc}
\usepackage[latin9]{inputenc}
\usepackage{geometry}
\usepackage{color}
\usepackage{babel}
\usepackage{amssymb}
\usepackage{amsmath}
\usepackage{amsthm}
\usepackage[numbers]{natbib}
\usepackage[unicode=true,pdfusetitle,
 bookmarks=true,bookmarksnumbered=false,bookmarksopen=false,
 breaklinks=false,pdfborder={0 0 1},backref=false,colorlinks=true]
 {hyperref}
\hypersetup{
 linkcolor=blue, urlcolor=marineblue, citecolor=blue, pdfstartview={FitH}, hyperfootnotes=false, unicode=true}

\makeatletter
\numberwithin{equation}{section}
\numberwithin{figure}{section}
\theoremstyle{plain}
\newtheorem{thm}{\protect\theoremname}
\theoremstyle{plain}

\newtheorem{prop}[thm]{\protect\propositionname}

\newcommand{\CC}{\mathbb{C}}
\newcommand{\Cc}{\cal{P}}
\newcommand{\Pp}{\cal{C}}
\usepackage{lettrine} 
\let\myFoot\footnote
\renewcommand{\footnote}[1]{\myFoot{#1\vspace{3mm}}}
\usepackage{lmodern}
\usepackage[T1]{fontenc}
\usepackage{concmath}
\usepackage{ae,aecompl}
\usepackage[T1]{fontenc}
\usepackage{amsmath}
\usepackage{graphicx}

\usepackage{tikz}

\makeatother

\providecommand{\corollaryname}{Corollary}
\providecommand{\propositionname}{Proposition}
\providecommand{\theoremname}{Theorem}

\begin{document}
\title{Remarks on results by M\"uger and Tuset on the moments of polynomials}

\author{Greg Markowsky\\
{\it Department of Mathematics,}\\
{\it Monash University, Melbourne, Australia}\\greg.markowsky@monash.edu.au \\
\\
Dylan Phung\\
{\it Yarra Valley Grammar School, Ringwood, Australia}\\
dkp1@yvg.vic.edu.au\\
}

\maketitle

\begin{abstract}
Let $f(x)$ be a non-zero polynomial with complex coefficients, and $M_p = \int_{0}^1 f(x)^p dx$ for $p$ a positive integer. In a recent paper, M\"uger and Tuset showed that $\limsup_{p \to \infty} |M_p|^{1/p} > 0$, and conjectured that this limit is equal to the maximum amongst the critical values of $f$ together with the values $|f(0)|$ and $|f(1)|$. We give an example that shows that this conjecture is false. It also may be natural to guess that $\limsup_{p \to \infty} |M_p|^{1/p}$ is equal to the maximum of $|f(x)|$ on $[0,1]$. However, we give a counterexample to this as well. We also provide a few more guesses as to the behaviour of the quantity $\limsup_{p \to \infty} |M_p|^{1/p}$.
\end{abstract}


\vspace{12pt}
Given a complex-valued function $f$ defined on the interval $[0,1]$, let the $p$-th moment of $f$ be defined as $M_p = \int_{0}^1 f(x)^p dx$ for $p$ a positive integer. In the elegant recent paper \cite{muger2020moments}, the following result was proved.

\begin{thm} \label{known}
Suppose $f$ is a polynomial with complex coefficients that is not identically 0. Then $L(f):=\limsup_{p \to \infty}\mid M_{p}\mid^{1/p} > 0$.
\end{thm}

Theorem \ref{known} has since been generalized to real analytic functions in \cite{estrada2020moments}. Within the proof in \cite{muger2020moments}, the authors defined the set ${\mathcal S}$ for polynomial $f$ as follows:

$$
{\mathcal S} = \{f(z) | z \in \CC, f'(z) = 0\} \cup \{f(0), f(1)\}.
$$

They then conjectured that perhaps

$$
\limsup_{p \to \infty}\mid M_{p}\mid^{1/p} = \max\{|s|: s \in {\mathcal S}\}.
$$

However, this conjecture is false, as we now show.

\begin{prop} 
Let $f(x) = 4 - (x+1)^2.$ Then $L(f) \leq 3$, but $\max\{|s|: s \in {\mathcal S}\} = 4$.
\end{prop}

{\bf Proof:} Clearly $|f(x)| \leq 3$ on $[0,1]$, and thus $L(f) \leq 3$ (in fact, $L(f) = 3$, since $\sup\{|f(x)|: x \in [0,1]\}=3$ and $f$ has real coefficients, by Proposition 1 in \cite{muger2020moments}). However, $f'(-1) = 0$, and $f(-1) = 4$. \qedsymbol

\vspace{8pt}

The authors of \cite{muger2020moments} also noted the simple fact that when the coefficients of $f$ are real we have $L(f) = \sup_{x \in [0,1]} |f(x)|$. It seems therefore natural to question whether the relation $L(f) =\sup_{x \in [0,1]} |f(x)|$ persists when $f$ is allowed to have complex coefficients. We will now give an example that shows that this is not true. We note first that we may change the interval in question to any other interval without changing the validity of the conjecture, since we may pass back and forth between them by linear maps. We will change it to $[-1,1]$, as the added symmetry will rather help our arguments, and make the natural changes to the definitions of $M_p$ and $L(f)$.  As noted, any counterexample to the statement in question must have at least one non-real coefficient, and our counterexample is the quadratic $f(x) = 1-(x+\frac{i}{2})^2$. There are doubtless many other examples, and we make no claims as to the optimality of our choice, but its simple form does allow us to analyze it easily. To be precise, we will prove the following.

\begin{prop} 
$\sup_{x \in [-1,1]} |f(x)| = |f(0)| = \frac{5}{4}$, but $L(f) \leq \frac{\sqrt{17}}{4}$.
\end{prop}

As a first step, let us understand the image $f([-1,1])$ in $\CC$, as well as the relevant branch of $f^{-1}$. The map $f_1(z) = (z+\frac{i}{2})^2$ maps the interval $[-1,1]$ to the sub-arc of the parabola $x=y^2-\frac{1}{4}$ between the points $\frac{3}{4}+i$ and $\frac{3}{4}-i$. To define $f_1^{-1}(w) = \sqrt{w} - \frac{i}{2}$, it is best to choose the branch cut for the square root to lie on the non-negative real line, with a point $re^{i \theta}$ with $0< \theta < 2 \pi$ being mapped to $r^{1/2}e^{i \theta/2} - \frac{i}{2}$. The image of $[-1,1]$ under $f(z) = 1-f_1(z)$ is therefore the reflection of this parabolic arc about the point $\frac{1}{2}$, and this is the sub-arc of the parabola $x=-y^2 + \frac{5}{4}$ connecting $\frac{1}{4}-i$ and $\frac{1}{4}+i$. For the inverse $f^{-1}(w) = \sqrt{1-w} - \frac{i}{2}$, the branch cut has also been reflected about the point $\frac{1}{2}$, and is the half-line $(-\infty, 1]$.

Let us now consider the moment $M_p = \int_{-1}^1 f(z)^p dz$. We change variables by $w = f(z)$, so that $dw = f'(z)dz$, whence

$$
dz = \frac{dw}{f'(f^{-1}(w))} = \frac{dw}{-2\sqrt{1-w}},
$$

where $\sqrt{\cdot}$ denotes the branch of the square root described above. We therefore obtain

$$
M_p = \int_{\Cc} \frac{w^p}{-2\sqrt{1-w}} dw,
$$

where $\Cc$ is the sub-arc of the parabola $x=-y^2 + \frac{5}{4}$ connecting $\frac{1}{4}-i$ and $\frac{1}{4}+i$. Cauchy's Theorem allows us to continuously deform $\Cc$ however we like, provided that the endpoints remain fixed and we do not pull it across a singularity or branch cut of the integrand. We may therefore replace $\Cc$ by $\Pp$, which is the sub-arc of the circle $\{x^2 + y^2 = \frac{17}{16}\}$ connecting $\frac{1}{4}-i$ and $\frac{1}{4}+i$. 

\begin{figure}
\centering{}~~~~~~~~~~~~~~\includegraphics[width=7cm,height=7cm,keepaspectratio]{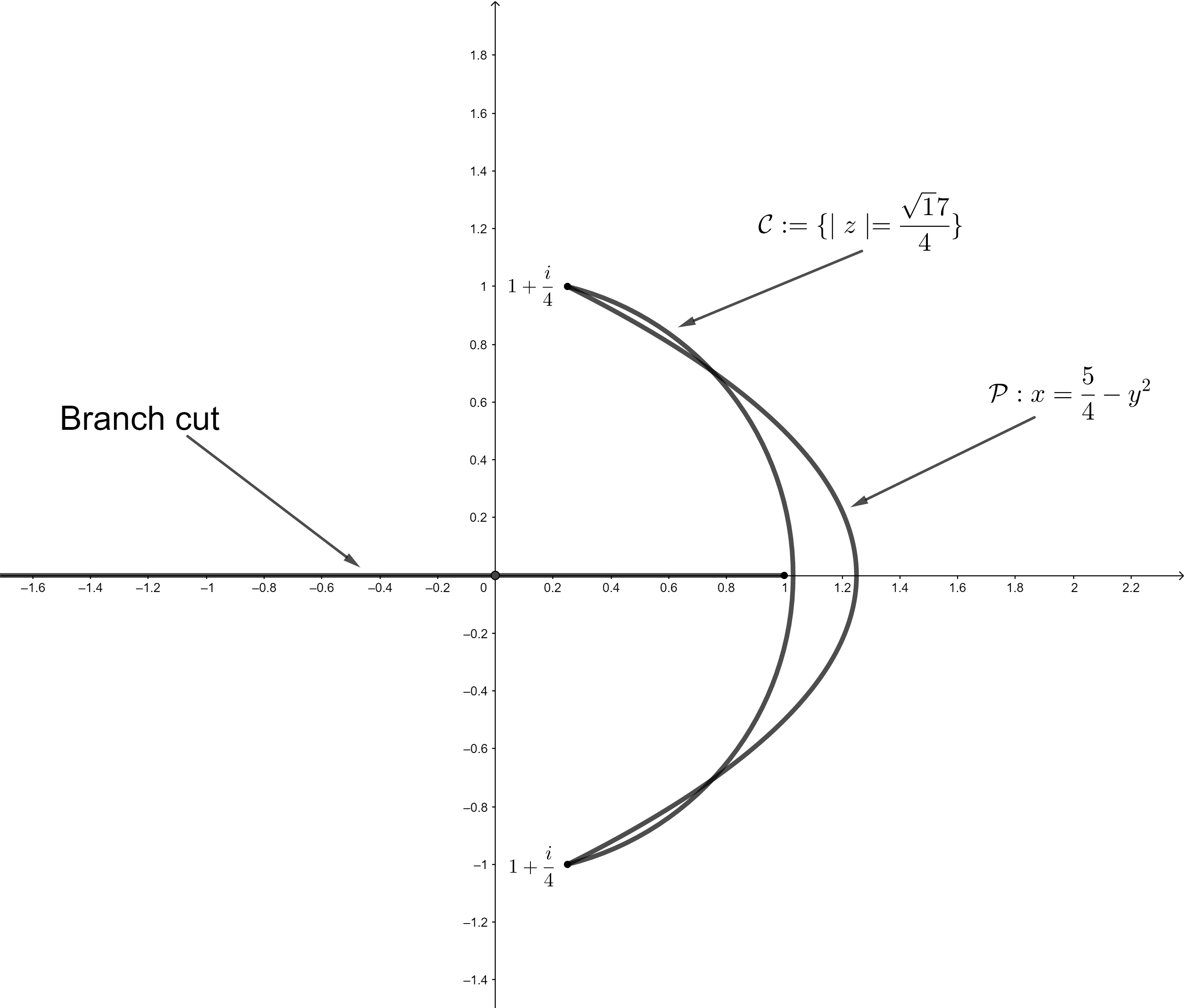}%\caption{The curves $\Pp$ and $\Cc$, and the branch cut. }
\end{figure}

We may now apply the standard $ML$-inequality (see \cite[Sec. 4.4]{schaum} or virtually any other complex analysis text) to obtain

$$
|M_p| \leq \mbox{length$(\Pp$)} \sup_{w \in \Pp} \Big|\frac{w^p}{-2\sqrt{1-w}}\Big| \leq K \frac{17^{p/2}}{4^p},
$$

where $K$ is a positive finite constant (note that $\frac{1}{2|\sqrt{1-w}|}$ is uniformly bounded and $|w| = \frac{\sqrt{17}}{4}$ on $\Pp$). Taking the $p$-th root and letting $p \to \infty$ yields the proposition. \qedsymbol

We remark that it seems to us that this general method, expressing the moments of $f$ as an integral over $f([-1,1])$ and then deforming the contour, would serve to create many more counterexamples to the conjecture. However, several potentially interesting new questions are also now indicated. For any polynomial $g$ on $[-1,1]$, let $\kappa(g) = \frac{L(g)}{\sup_{x \in [-1,1]} |g(x)|}$. It is clear that $\kappa(g) \leq 1$, and we have given an example $f$ with $\kappa(f) \leq \frac{\sqrt{17}}{5}$.

\vspace{8pt}

{\bf Question 1:} Is there a universal constant $\delta>0$, such that $\kappa(g)>\delta$ for all non-zero polynomials $g$?

\vspace{8pt}

{\bf Question 2:} Failing that, is there a constant $\delta(n)>0$, depending on $n$, such that $\kappa(g)>\delta(n)$ for all non-zero polynomials $g$ of degree $n$?

\vspace{8pt}

{\bf Question 3:} What are the optimal constants in either of the above questions?

\vspace{8pt}

{\bf Question 4:} Is it possible to find a polynomial $f$ such that $L(f) < \max\{|f(a)|,|f(b)|\}$, where $[a,b]$ is the interval that the moments are taken over? 

\vspace{8pt}

For this final question, we remark that our method does not seem to apply here, since we are not able to move the endpoints. If such a polynomial exists then a new method for proving it will need to be devised. We should also remark that a number of possibilities for such a polynomial have already been eliminated by the comments at the end of \cite{muger2020moments}. 

%\section{Acknowledgements}

%We would like to thank Maher Boudabra for useful conversations, and for helping with the image.

\bibliographystyle{plain}
\bibliography{biblio}

\end{document}